\documentclass[10pt,aps,prb,twocolumn, nofootinbib]{revtex4-1}
\usepackage{mathtools,amsmath}
\usepackage{amssymb,amsmath}
\usepackage{subfig}
\usepackage{graphicx}
\usepackage{color}
\usepackage{hyperref}
\usepackage{listings}

\hypersetup{
    bookmarks=true,         
    unicode=false,          
    pdftoolbar=true,        
    pdfmenubar=true,        
    pdffitwindow=false,     
    pdfstartview={FitH},    
    pdfauthor={Kevin Multani},     
    colorlinks=true,       
    linkcolor=blue,          
    citecolor=red,        
}
\graphicspath{{figures/}}

\begin{document}

\definecolor{dkgreen}{rgb}{0,0.6,0}
\definecolor{gray}{rgb}{0.5,0.5,0.5}
\definecolor{mauve}{rgb}{0.58,0,0.82}

\lstset{frame=tb,
  	language=Matlab,
  	aboveskip=3mm,
  	belowskip=3mm,
  	showstringspaces=false,
  	columns=flexible,
  	basicstyle={\small\ttfamily},
  	numbers=none,
  	numberstyle=\tiny\color{gray},
 	keywordstyle=\color{blue},
	commentstyle=\color{dkgreen},
  	stringstyle=\color{mauve},
  	breaklines=true,
  	breakatwhitespace=true
  	tabsize=3
}

\title{The half--form $\sqrt{dx}$}
\author{S. Camosso}
\affiliation{Department of Mathematics and Phisics of \\ ``A.Bertoni--G.Soleri'' Secondary School}
\date{\today}

\begin{abstract}
In this brief article we try to find  an ``interpretation'' for  the formalism $\sqrt{dx}$.

\end{abstract}

\maketitle

The problem to make sense of $\int \sqrt{dx}$ is treated in this article. In what follows we explore a solution called ``the corrected'' integral, formally we denote it by $\prescript{}{\gamma}{\int_{[a,b]}}\sqrt{dx}$. Some inspiration for this study can be suggested by Ramanujan \cite{uno} from its convergent series. The problem is approached from the point of view of the calculus, with some relations to the fractional calculus in Ross \cite{otto}. 
A well definition of the symbol $\sqrt{dx}$ can be useful in geometric quantization (Woodhouse \cite{quattro}) and different areas of mathematics. In the contest of geometric quantization the symbol $\sqrt{dx}$ denotes a section of the quantum line bundle defined on a compact, complex, symplectic manifold $M$. 

\section{\label{sec:one}Integration of $\sqrt{dx}$ with a Riemannian integral}
In this section we start to examine the situation evaluating the following integral:

\begin{equation}
\label{sqrt{dx}1}
\int_{[a,b]}\,\sqrt{dx}.
\end{equation}

Let $\mathcal{P}\,=\,\{x_{0}=a,x_{1},\ldots,x_{n}=b\}$ be a partition of the interval $[a,b]$ in $n$ subintervals of amplitude $\frac{b-a}{n}$. Thus

\begin{equation}
\label{sqrt{dx}2}
\int_{[a,b]}\sqrt{dx}\,=\,\lim_{n\rightarrow +\infty} \sum_{i=1}^{n}\sqrt{\Delta x_{i}},
\end{equation}
\noindent 
where $\Delta x_{i}\,=\, x_{i}-x_{i-1}$ for $i=1, \ldots, n$. This is the Riemann integral where, instead to consider as ``base of rectangles'' the quantities $\Delta x_{i}$, we consider $\sqrt{\Delta x_{i}}$. The result is:

\begin{equation}
\label{sqrt{dx}3}
\begin{multlined}
\int_{[a,b]}\,\sqrt{dx}\,=\,\lim_{n\rightarrow +\infty} \sum_{i=1}^{n}\sqrt{\frac{b-a}{n}}\,=\, \\ 
\,=\, \lim_{n\rightarrow +\infty}\sqrt{(b-a)n}\,=\,+\infty.
\end{multlined}
\end{equation}

The integral diverges at $+\infty$, the result is not satisfactory at all. 

\section{\label{sec:two}On a Ramanujan sum}

In order to find a way to make this infinity ``disappear'', let us consider this result due to Ramanujan\cite{uno}:

\begin{equation}
\label{Ramanujan1}
\begin{multlined}
\sqrt{1}+\sqrt{2}+\sqrt{3}+\cdots +\sqrt{n}\,=\, \\ \,=\, \frac{2}{3}n\sqrt{n}+\frac{1}{2}\sqrt{n}-\frac{\zeta\left(\frac{3}{2}\right)}{4\pi}+\frac{1}{24\sqrt{n}}+ \mathit{o}\left(\frac{1}{n^{\frac{5}{2}}}\right).
\end{multlined}
\end{equation}

This is an asymptotic expansion of the sum of the square roots of the first $n$ natural numbers. The main term in the expansion, when $n$ goes to infinity, is $\frac{2}{3}n^{\frac{3}{2}}$. Another similar series, always due to Ramanujan is:

\begin{equation}
\label{Ramanujan2}
\begin{multlined}
\frac{1}{\sqrt{1}}+\frac{1}{\sqrt{2}}+\frac{1}{\sqrt{3}}+\cdots +\frac{1}{\sqrt{n}}\,=\, \\ \,=\, 2\sqrt{n}+\frac{1}{2\sqrt{n}}+ \zeta\left(\frac{1}{2}\right)+\mathit{o}\left(\frac{1}{n^{\frac{3}{2}}}\right).
\end{multlined}
\end{equation}

The idea can be to modify the sum of rectangles using the factor $\gamma(n)\,=\,\frac{1}{\sqrt{n}}$. In this case:

\begin{equation}
\label{solution_sqrt{dx}3}
\begin{multlined}
\int_{[a,b]}\,\sqrt{dx}\,=\,\lim_{n\rightarrow +\infty} \gamma(n)\cdot \sum_{i=1}^{n}\sqrt{\frac{b-a}{n}}\,=\, \\ 
\,=\, \sqrt{(b-a)}\lim_{n\rightarrow +\infty}\frac{\sqrt{n}}{\sqrt{n}}\,=\, \sqrt{b-a}.
\end{multlined}
\end{equation}

This is not the Riemannian integral used before but it is a variation with a normalized sum. The problem was in fact the divergence of the series of square roots. This sort of ``correction'' bring to a finite result that corresponds to the square root of the initial interval. 

We denote this corrected integral with the notation: 

\begin{equation}
\prescript{}{\gamma}{\int_{[a,b]}}\,\sqrt{dx}\,=\,\lim_{n\rightarrow +\infty} \prescript{}{\gamma}{\sum}_{i=1}^{n}\sqrt{\Delta x_{i}},
\end{equation}
\noindent 
where $\prescript{}{\gamma}{\sum}$ is the sum corrected by the factor $\gamma(n)$ and $\prescript{}{\gamma}{\int}$ is the ``corrected'' integral.

\section{\label{sec:three}The geometric quantization program and the definition of $\sqrt{dx}$}

An attempt to define the half--form $\sqrt{dx}$ has been done in the program of geometric quantization. The geometric quantization is a process that associate to a symplectic manifold an Hilbert space representing the space of quantum states (for references about the geometric quantization program see Woodhouse \cite{quattro} or works of Kostant \cite{sette} and Souriau \cite{sei}). 

The idea is to see the half--form $\sqrt{dx}$ as a section of the square root of the canonical line bundle associated to the polarization adopted during the process of geometric quantization. The concept of $\frac{1}{2}$--form is strictly correlated to the concept of $\frac{1}{2}$--density. In particular we have that $\frac{1}{2}$--densities can be identified with $\frac{1}{2}$--forms. Let us consider the case of $\mathbb{R}$. The space $\mathbb{R}$ is a vectorial space where we can choose an orientation. If $F(\mathbb{R})$ represents the set of frames of $\mathbb{R}$ we can consider the action of $GL(1,\mathbb{R})$ on $F(\mathbb{R})$ that is simply the scalar multiplication. We define the set of $\frac{1}{2}$-densities as:

\begin{equation}
\label{densities}
\begin{multlined}
\left|\mathbb{R}\right|^{\frac{1}{2}}\,=\, \left\{ \nu:F(\mathbb{R})\rightarrow \mathbb{C}: \nu(a\cdot v)\,=\, \nu(v)|\det{a}|^{\frac{1}{2}},\right. \\ \left.\forall v\in F(\mathbb{R}),a\in GL(1,\mathbb{R})\right\}.
\end{multlined}
\end{equation}

This set is a line bundle denoted by $\left|\mathbb{R}\right|^{\frac{1}{2}}\rightarrow \mathbb{R}$ where the $1$--dimensional fiber at $x\in\mathbb{R}$ is $\left|\mathbb{R}_{x}\right|^{\frac{1}{2}}\rightarrow \mathbb{R}_{x}$. 
Rigorous definitions of $\frac{1}{2}$--forms and $\frac{1}{2}$--densities with its properties can be found in Guillemin and Sternberg \cite{tre}, Hall \cite{due} and Rawnsley \cite{cinque}.

In geometric quantization the $\frac{1}{2}$--forms define, intrinsically, an half--form Hilbert space with an inner product and a norm. Let us consider the case of $\mathbb{R}$. Then the symplectic manifold is $M\,=\,T^{\vee}\mathbb{R}\,=\,\mathbb{R}^{2}$. Let $P$ be the vertical polarization of $M$ with the orientation of $\mathbb{R}$ so that oriented $1$--forms are positive multiple of $dx$. Let $\sqrt{K}_{P}$ to be the trivial bundle with trivializing section $\sqrt{dx}$ such that $\sqrt{dx}\otimes \sqrt{dx}\,=\, dx$. Then the half--form $\sigma\,=\, f(x)\sqrt{dx}$, for some real function $f(x)$, has the norm:

\begin{equation}
\label{norm}
\|\sigma\|^2\,=\, \int_{\mathbb{R}}\left|f(x)\right|^2 dx.
\end{equation}

\section{\label{sec:three} The corrected integral as application from the $\frac{1}{2}$--forms to $\mathbb{R}$}

Let us consider the following $\frac{1}{2}$--form $\sigma\,=\,\sqrt{dx}$. In order to be precise, we must view this form as the form $1\otimes \sqrt{dx}$ of the quantum line bundle $L\otimes \sqrt{K}_{P}$. In this case we have that:

\begin{equation}
\label{norm_example}
\|\sigma\|^2\,=\, \int_{[a,b]} dx\,=\, b-a,
\end{equation}
\noindent 
where we considered as a base space the interval $[a,b]$. Now the question is if exists a square root of the following equation:

\begin{equation}
\label{norm_example_2}
\sigma\cdot\sigma\,=\, b-a,
\end{equation}
\noindent 
where the product here is the squared norm in the Hilbert space of half--forms. In order to answer the question let us consider the corrected integral:

\begin{equation}
\label{solution_problem}
\begin{multlined}
\prescript{}{\gamma}{\int_{[a,b]}}\,\sqrt{dx}\,=\,\sqrt{b-a}.
\end{multlined}
\end{equation}

We can see that:

\begin{equation}
\label{solution_problem_2}
\begin{multlined}
\prescript{}{\gamma}{\int_{[a,b]}}\,\sqrt{dx}\cdot  \prescript{}{\gamma}{\int_{[a,b]}}\,\sqrt{dx}\,=\, (b-a).
\end{multlined}
\end{equation}

So it is possible to apply the corrected integral in order to find the corrected result. In other terms we can see the corrected integral as a map $\prescript{}{\gamma}{\int}:\left|\mathbb{R}\right|^{\frac{1}{2}}\rightarrow \mathbb{R}$, omitting the quantum line bundle associated. 

\section{\label{sec:three} Interesting integrals in $\sqrt{dx}$}

Let us consider the integral:

\begin{equation}
\label{sqrt{dx}2x}
\int_{[a,b]}x\,\sqrt{dx}.
\end{equation}

Without the correction the integral diverges. We can try with:

\begin{equation}
\label{sqrt{dx}2xcorr}
\prescript{}{\gamma}{\int_{[a,b]}}x\,\sqrt{dx}.
\end{equation}

Let us consider the following partition of $[a,b]$ with $x_{i}\,=\, a+(i-1)\cdot h$ and $x_{i+1},=\, a+i\cdot h$, where $h=\frac{b-a}{n}$. Thus:

\begin{equation}
\label{sqrt{dx}2xcorr3}
\begin{multlined}
\prescript{}{\gamma}{\int_{[a,b]}}x\,\sqrt{dx}\,=\, \lim_{n\rightarrow +\infty}\gamma(n)\sum_{i=1}^{n}\left(a+i\cdot h\right)\sqrt{h}\,=\, \\ \,=\, \lim_{n\rightarrow +\infty} a\cdot \sqrt{b-a} + \frac{\sqrt{(b-a)^3}}{n^{2}}\sum_{i=1}^{n}i\,=\, \\ \,=\, a\cdot \sqrt{b-a}+\frac{\sqrt{(b-a)^3}}{2},
\end{multlined}
\end{equation}
\noindent 
where we used the fact that $\sum_{i=1}^{n}i\,=\, \frac{n\cdot (n+1)}{2}$. 

A similar result is the following:

\begin{equation}
\label{sqrt{dx}2xquadrato}
\begin{multlined}
\prescript{}{\gamma}{\int_{[a,b]}}x^2\,\sqrt{dx}  \,=\, a^2\sqrt{b-a}+\\+(b-a)\sqrt{b-a}+\frac{1}{3}(b-a)^2\sqrt{b-a}.
\end{multlined}
\end{equation}

The calculations are similar to the previous case where the formula used now is $\sum_{i=1}^{n}i^2\,=\, \frac{n\cdot (n+1)\cdot (2n+1)}{6}$.

\section{\label{sec:four}  $\sqrt{dx}$ from Fractional Calculus}

In this section we see as our definition of corrected fractional integral is in perfect agreement with the definition from the fractional calculus \cite{otto}. We start recalling the formula for the $\frac{1}{2}$--integral between $a$ and $b$, this is given by the formula:

\begin{equation}
\label{fractional_calculus1}
D^{-\frac{1}{2}}_{[a,b]}f(t)\,=\, \frac{1}{\Gamma\left(\frac{1}{2}\right)}\int_{a}^{b}(b-t)^{-\frac{1}{2}}f(t)dt
\end{equation}

Now we observe that if we compute the $\frac{1}{2}$--integral for the constant function $f(t)=1$ we find that:

\begin{equation}
\label{fractional_calculus2}
D^{-\frac{1}{2}}_{[a,b]}1\,=\, \frac{2}{\Gamma\left(\frac{1}{2}\right)}\sqrt{b-a}.
\end{equation}

We observe that we have considered the fractional integral over the interval $[a,b]$. Usually the integral $(\ref{fractional_calculus1})$ is considered on the intervall $[0,x]$ and, for general calculations of $(\ref{fractional_calculus1})$, we need to use the beta integral:

\begin{equation}
\label{fractional_calculus3}
\int_{0}^{x}(x-y)^{d}y^bdy\,=\, \frac{\Gamma(d+1)\Gamma(b+1)}{\Gamma(b+d+2)}x^{b+d+1}.
\end{equation}

We can compare the result from the theory of fractional calculus $\frac{2}{\Gamma\left(\frac{1}{2}\right)}\sqrt{b-a}$ with our definition using the corrected integral $(\ref{solution_problem})$ that gives $\sqrt{b-a}$ (only a constant factor of difference!). In fact:

\begin{equation}
\label{fractional_calculus2mod}
\prescript{}{\gamma}{\int_{[a,b]}}1\,\sqrt{dx} \,=\, \frac{\Gamma\left(\frac{1}{2}\right)}{2}D^{-\frac{1}{2}}_{[a,b]}1.
\end{equation}

\section{\label{sec:five}  Observations and conclusion}

The main observation is that we have realized the following relation:

\begin{equation}
\label{intercange}
\prescript{}{\gamma}{\int_{[a,b]}}1\,\sqrt{dx} \,=\, \sqrt{\int_{[a,b]}1\,dx}.
\end{equation}

The relation seems not true in general $\prescript{}{\gamma}{\int_{[a,b]}}f(x)\,\sqrt{dx} \,\not=\, \sqrt{\int_{[a,b]}f(x)\,dx}$. The fact is clear observing that for $f(x)=x$ the previous relation is false.
We can try to define an integral function $F(x)=\prescript{}{\gamma}{\int_{[0,x]}}1\,\sqrt{ds}$, in this case:

\begin{equation}
\label{intercange}
F(x)\,=\, \prescript{}{\gamma}{\int_{[0,x]}}1\,\sqrt{ds} \,=\, \sqrt{x}.
\end{equation}

Other questions are open, it is possible a similar definition for $\sqrt[n]{dx}$ ? (for $n=3,4, \ldots $). Another important discussion theme concerns the geometrical meaning of this correction.

\end{document}